\RequirePackage{fix-cm}
\documentclass[smallextended]{svjour3}       
\smartqed  
\usepackage{graphicx}
\usepackage{amsmath,amsfonts,bm,mathtools}
\begin{document}

\title{Global solutions to a class of CEC benchmark constrained optimization problems
}


\author{Xiaojun Zhou         \and
        David Yang Gao \and
        Chunhua Yang 
}


\institute{Xiaojun Zhou \at
              School of Science, Information Technology and Engineering, University of Ballarat, Victoria 3353, Australia.\\
              School of Information Science and Engineering, Central South University, Changsha 410083, China.\\
           \and
           David Yang Gao \at
              School of Science, Information Technology and Engineering, University of Ballarat, Victoria 3353, Australia.\\
           \and
           Chunhua Yang \at
              School of Information Science and Engineering, Central South University, Changsha 410083, China.
}

\date{Received: date / Accepted: date}

\maketitle

\begin{abstract}
This paper aims to solve a class of CEC benchmark constrained optimization problems that have been widely studied by nature-inspired optimization algorithms. Global optimality condition based on canonical duality theory is derived. Integrating the dual solutions with the KKT conditions, we are able to obtain the approximate solutions or global solutions easily.

\end{abstract}

\section{Introduction}
Nature-inspired optimization algorithms, such as genetic algorithm (GA), evolution strategy (ES), particle swarm optimization (PSO) and differential evolution (DE), have received considerable attention in recent decades due to their strong adaptability and easy implementation. Strictly speaking, these algorithms are unconstrained optimization procedures, and therefore it is necessary to
find techniques to deal with the constraints when solving constrained optimization problems. The most common approach to handle constraints is the penalty function method.
The idea of this approach is to transform a
constrained optimization problem into an unconstrained one by adding a certain term
to the objective function based on the amount of constraint violation.
Then, some special representations and operators are designed to preserve the feasibility of solutions at all times or to repair a solution when it is infeasible. Multiobjective optimization techniques are also applied to manage constraints. The main idea is to rewrite the single objective optimization problem as a multiobjective optimization problem in which the constraints in original problem are treated as additional objectives \cite{Coello,Michalewicz}.\\
\indent By introducing a Lagrangian multiplier vector to relax the constraints, the classical Lagrangian is a saddle function if the objective function and all of the constraints are convex. Under certain constraint qualifications, the well-known strong \textit{min-max} duality relation holds, and in this case, the problem can be easily solved by well-developed convex programming techniques. However, due to nonconvexity of either objective function or the constraints, the Lagrangian is no longer a saddle function and only the weak duality relation holds, leading to the \textit{duality gap} in global optimization \cite{Bazaraa,Boyd,GRS2009}. In order to bridge the gap inherent in the classical Lagrange duality theory, the \textit{canonical duality theory} has been developed recently. Its core is to transform a nonconvex primal minimization problem to the concave canonical dual maximization problem over a convex space without duality gap by a canonical dual transformation \cite{Gao2009}.\\
\indent As shown by the global optimality condition contained in the canonical duality theory, if the dual solution is in the positive definite domain, it is easy to get the corresponding global solution to the primal problem. However, if the condition is not satisfied, some strategies are necessary to recover the global solution \cite{Fang-gaoetal07,Gao2010,Wangetal}. In this paper, we focus on solving a class of CEC (Congress on Evolutionary Computation) benchmark constrained optimization problems that have been widely studied by nature-inspired algorithms. By integrating the canonical duality theory with the KKT conditions, we are able to obtain the approximate solutions or global solutions easily.

\section{The canonical duality theory}
In this paper, we focus on the following quadratic optimization problem with quadratic and box constraints (primal problem):
\begin{eqnarray}\label{eq1}
(\mathcal{P}): \min && \Big\{ P(\mathbf{x})  = \frac{1}{2}\mathbf{x}^T A \mathbf{x} - \mathbf{a}^T\mathbf{x} - a: \mathbf{x} \in \mathbb{R}^n \Big\},   \nonumber \\
\mathrm{s.t.} && \mathbf{g}(\mathbf{x}) = \{g_j(\mathbf{x})\} = \{\frac{1}{2}\mathbf{x}^T B_j \mathbf{x} - \mathbf{b}_j^T\mathbf{x} - b_j\} \leq \mathbf{0}, j = 1, \cdots, m,  \nonumber \\
&&  c_i \leq x_i \leq d_i, i = 1, \cdots, n,
\end{eqnarray}
where, $\mathbf{x} = (x_1, \cdots, x_n), A = A^T, B_j = B_j^T \in \mathbb{R}^{n \times n}$ are symmetric matrices, $\mathbf{a}, \mathbf{b}_j \in \mathbb{R}^n$ are given vectors, $a,b_j,c_i, d_i$ are constant.\\
\indent Let $E_i \in \mathbb{R}^{n \times n}$, $\mathbf{e}_i \in \mathbb{R}^{n}$ be a diagonal matrix and a unit vector, with all zeros except a one in the position $(i,i)$ and $(i)$, respectively, $B_k = 2E_k$, $\mathbf{b}_k = (c_k+d_k)\mathbf{e}_k$, $b_k = c_kd_k$, $k = m+1, \cdots, m+n$,
then constraints in $(\mathcal{P})$ can be uniformly rewritten to
\begin{eqnarray}\label{eq2}
\mathbf{g}(\mathbf{x}) = \{g_k(\mathbf{x})\} = \{\frac{1}{2}\mathbf{x}^T B_k \mathbf{x} - \mathbf{b}_k^T\mathbf{x} - b_k\} \leq \mathbf{0}, k = 1, \cdots, m+n.
\end{eqnarray}
\indent Firstly, we introduce an indicator function:
\begin{equation} \label{eq3}
W(\bm{\epsilon}) =
\left\{ \begin{aligned}
 0       & ~~\mathrm{if} ~\bm{\epsilon} \leq \mathbf{0} \\
 +\infty & ~~\mathrm{otherwise}
\end{aligned} \right.
\end{equation}
where, $\bm{\epsilon} = (\epsilon_1, \cdots, \epsilon_{m+n})$, and let $U(\mathbf{x}) = -f(\mathbf{x})  = - \frac{1}{2}\mathbf{x}^T A \mathbf{x} + \mathbf{a}^T\mathbf{x} + a$, then the primal problem ($\mathcal{P}$) can be written in the following form:
\begin{equation}\label{eq4}
    \min \{\mathrm{\Pi}_0(\mathbf{x}) = W(\mathbf{g}(\mathbf{x})) - U(\mathbf{x}): \mathbf{x} \in \mathbb{R}^n\}.
\end{equation}
\indent Secondly, we introduce a nonlinear operator
\begin{equation}\label{eq5}
    \bm\xi = (\xi_1, \cdots, \xi_{m+n}) = \{\xi_k\} = \mathrm{\Lambda}(\mathbf{x}) = \{\frac{1}{2}\mathbf{x}^T B_k \mathbf{x} - \mathbf{b}_k^T\mathbf{x} - b_k\}
\end{equation}
so that $\mathbf{g}(\mathbf{x})$ can be recast by $V(\bm\xi) = \bm\xi = \mathbf{g}(\mathbf{x})$, and then the primal problem ($\mathcal{P}$) can be reformed in the canonical form:
\begin{equation}\label{eq6}
    \min \{\mathrm{\Pi}(\mathbf{x}) = W(V(\mathrm{\Lambda}(\mathbf{x}))) - U(\mathbf{x}): \mathbf{x} \in \mathbb{R}^n\}.
\end{equation}
\indent By the Fenchel transformation, the conjugate function $W^{\sharp}(\bm{\sigma})$ of $W(\bm{\epsilon})$ can be defined by
\begin{equation}\label{eq7}
   W^{\sharp}(\bm{\sigma}) = \underset{\bm{\epsilon}}{\mathrm{sup}}\{\bm{\epsilon}^T \bm{\sigma} - W(\bm{\epsilon}) \} =
\left\{ \begin{aligned}
 0      & ~~\mathrm{if} ~ \bm{\sigma} \geq 0\\
 +\infty & ~~\mathrm{otherwise}
\end{aligned} \right.
\end{equation}
which is convex and lower semi-continuous. \\
\indent According to the relations $\bm{\sigma} \in \partial W(\bm{\epsilon})\Leftrightarrow \bm{\epsilon} \in \partial W^{\sharp}(\bm{\sigma}) \Leftrightarrow W(\bm{\epsilon}) + W^{\sharp}(\bm{\sigma}) = \bm{\epsilon}^T\bm{\sigma}$ from convex analysis, we can replace $W(\mathbf{g}(\mathbf{x}))$ by $\mathbf{g}^T(\mathbf{x})\bm\sigma -  W^{\sharp}(\bm{\sigma})$, and then
we can get the extended Lagrangian
\begin{equation}\label{eq8}
    \Xi_0(\mathbf{x},\bm{\sigma}) = \mathbf{g}^T(\mathbf{x})\bm\sigma -  W^{\sharp}(\bm{\sigma}) - U(\mathbf{x}).
\end{equation}
\indent Next, we introduce the invertible duality mapping
\begin{equation}\label{eq9}
    \bm\varsigma = (\varsigma_1, \cdots, \varsigma_{m+n}) = \{\varsigma_k\} = \nabla V(\bm\xi) = I.
\end{equation}
\indent Defining the Legendre conjugate $V^{\ast}(\bm\varsigma) = \mathrm{sta}\{\bm\xi^T \bm\varsigma  - V(\bm\xi)\}$, and using the equivalent relations $\bm\varsigma  = \nabla V(\bm\xi)\Leftrightarrow \bm\xi = \nabla V^{\ast}(\bm\varsigma) \Leftrightarrow \bm\xi^T \bm\varsigma = V(\bm\xi) + V^{\ast}(\bm\varsigma) $,
we can replace $\mathbf{g}(\mathbf{x})$ by $\mathrm{\Lambda}^T(\mathbf{x}) \bm\varsigma - V^{\ast}(\bm\varsigma)$,
so we obtain the generalized complementary function
\begin{equation}\label{eq10}
   \Xi(\mathbf{x},\bm{\sigma}) = \frac{1}{2}\mathbf{x}^T G (\bm\sigma) \mathbf{x} - \mathbf{x}^T F(\bm\sigma) - \bm\sigma^T \mathbf{d} - a,
\end{equation}
where,
\begin{eqnarray*}
G (\bm\sigma) = A + \sum\limits_{k=1}^{m+n}\sigma_k B_k, F(\bm\sigma) = \mathbf{a} + \sum\limits_{k=1}^{m+n}\sigma_k \mathbf{b}_k,
\end{eqnarray*}
and
\begin{eqnarray*}
\mathbf{d} = (b_1,\cdots, b_{m+n})^T, \bm{\sigma} = (\sigma_1, \cdots, \sigma_{m+n}) \in \mathbb{R}^{m+n}_+ = \{\bm{\sigma} \in \mathbb{R}^{m+n}| \bm{\sigma} \geq 0 \}.
\end{eqnarray*}
\indent By using the generalized complementary function, the canonical dual function $P^d(\bm{\sigma})$ can be formulated as
\begin{eqnarray}
P^d(\bm{\sigma}) = \underset{\mathbf{x}}{\mathrm{sta}}\{\Xi(\mathbf{x},\bm{\sigma})\}.
\end{eqnarray}
\indent Solving the critical points of $\Xi(\mathbf{x},\bm{\sigma})$, we can get the canonical equilibrium equation
\begin{eqnarray}
G(\bm\sigma) \mathbf{x} = F(\bm\sigma).
\end{eqnarray}
\indent For any given $\bm\sigma$, if $F(\bm{\sigma})$ is in the column space of $G(\bm{\sigma})$, denoted by $\mathcal{C}_{ol}(G(\bm{\sigma}))$, i.e., a linear space spanned by the columns of $G(\bm{\sigma})$, the solution of the canonical equilibrium equation can be well defined by
\begin{eqnarray}
\mathbf{x} = G^{\dagger}(\bm{\sigma})F(\bm{\sigma}),
\end{eqnarray}
where, $G^{\dagger}(\bm{\sigma})$ denotes the Moore-Penrose generalized inverse of $G(\bm{\sigma})$.\\
\indent Then, the canonical dual function can be written explicitly as follows
\begin{eqnarray}
P^d(\bm{\sigma}) = -\frac{1}{2}F^T(\bm{\sigma})G^{\dagger}(\bm{\sigma})F(\bm{\sigma}) - \bm\sigma^T \mathbf{d} - a.
\end{eqnarray}
\indent Finally, the canonical dual problem can be posed as follows:
\begin{eqnarray}\label{pd}
(\mathcal{P}^d): \max \{P^d(\bm{\sigma}):\bm\sigma \in \mathcal{S}^+_a\}
\end{eqnarray}
where, the dual feasible space is defined by $\mathcal{S}^+_a = \{\bm\sigma \in \mathbb{R}^{m+n}_{+} | G(\bm{\sigma})\succeq 0\}$.\\

\noindent \textbf{Theorem 1} (\textbf{Global Optimality Condition}). Suppose $\bar{\bm\sigma}$ is a KKT point of ($\mathcal{P}^d$). If $G(\bar{\bm\sigma})\succ 0$, then $\bar{\mathbf{x}} = G^{-1}(\bar{\bm\sigma})F(\bar{\bm\sigma})$ is the global minimizer of ($\mathcal{P}$). If det$(G(\bar{\bm\sigma}))=0$, the global minimizer $\bar{\mathbf{x}}$ of ($\mathcal{P}$) is contained in the canonical equilibrium equation $G(\bar{\bm\sigma}) \bar{\mathbf{x}} = F(\bar{\bm\sigma})$.\\
\textbf{Proof}. By introducing Lagrange multiplier $\bm{\epsilon} \in \mathbb{R}_{-}^{m+n}$ (where $\mathbb{R}_{-}^{m+n}$ is the nonpositive orthant of $\mathbb{R}^{m+n}$) associated with $\bm{\sigma} \geq 0$, the Lagrangian $L(\bm{\epsilon}, \bm{\sigma})$ is given by
\begin{eqnarray}
L(\bm{\epsilon}, \bm{\sigma}) = -\frac{1}{2}F^T(\bm{\sigma})G^{\dagger}(\bm{\sigma})F(\bm{\sigma}) - \bm\sigma^T \mathbf{d} - a - \bm{\epsilon}^T \bm{\sigma}.
\end{eqnarray}
It is easy to prove that the criticality conditions $\nabla_{\bm{\sigma}} L(\bm{\epsilon}, \bm{\sigma}) = 0$ lead to
\begin{eqnarray}
\bm{\epsilon} =
\begin{pmatrix}
\epsilon_1\\
\cdots\\
\epsilon_{m+n}
\end{pmatrix}
=
\begin{pmatrix}
\frac{1}{2}\bar{\mathbf{x}}^T B_1 \bar{\mathbf{x}} - \mathbf{b}_1^T \bar{\mathbf{x}} -b_1 \\
\cdots\\
\frac{1}{2}\bar{\mathbf{x}}^T B_{m+n} \bar{\mathbf{x}} - \mathbf{b}_{m+n}^T \bar{\mathbf{x}} -b_{m+n}
\end{pmatrix}
\end{eqnarray}
and the accompanying KKT conditions include
\begin{eqnarray}
0 \leq \bar{\sigma}_k \perp \frac{1}{2}\bar{\mathbf{x}}^T B_k \bar{\mathbf{x}} - \mathbf{b}_k^T \bar{\mathbf{x}} -b_k \leq 0, k = 1, \cdots, m+n.
\end{eqnarray}
From the complementary slackness, we can see that the $\bar{\mathbf{x}}$ satisfies the constraints in ($\mathcal{P}$).
Furthermore, since $\bar{\bm{\sigma}} \geq 0$ for any $\mathbf{g}(\mathbf{x})\leq 0$, we have
\begin{eqnarray}
P(\mathbf{x}) &\geq&  P(\mathbf{x}) + \bar{\bm{\sigma}}^T \mathbf{g}(\mathbf{x}) \nonumber \\
&=&\frac{1}{2}\mathbf{x}^T A \mathbf{x} - \mathbf{a}^T\mathbf{x} - a + \sum_{k=1}^{m+n} (\frac{1}{2}\mathbf{x}^T \bar{\sigma}_k B_k  \mathbf{x} - \bar{\sigma}_k \mathbf{b}_k^T\mathbf{x} - \bar{\sigma}_k b_k) \nonumber \\
&=& \frac{1}{2}\mathbf{x}^T G (\bar{\bm{\sigma}}) \mathbf{x} - \mathbf{x}^T F(\bar{\bm{\sigma}}) - \bar{\bm{\sigma}}^T \mathbf{d} - a \nonumber \\
&=& \Xi(\mathbf{x},\bar{\bm{\sigma}}).
\end{eqnarray}
Noting that $P(\bar{\mathbf{x}}) =  \Xi(\bar{\mathbf{x}},\bar{\bm{\sigma}}), \nabla_{\mathbf{x}}\Xi(\bar{\mathbf{x}},\bar{\bm{\sigma}}) = 0$ and $\Xi(\mathbf{x},\bar{\bm{\sigma}})$ is a quadratic function with respect to $\mathbf{x}$, we have
\begin{eqnarray}
P(\mathbf{x}) - P(\bar{\mathbf{x}}) &\geq& \Xi(\mathbf{x},\bar{\bm{\sigma}}) - \Xi(\bar{\mathbf{x}},\bar{\bm{\sigma}}) \nonumber \\
&=& (\mathbf{x} - \bar{\mathbf{x}})\nabla_{\mathbf{x}}\Xi(\bar{\mathbf{x}},\bar{\bm{\sigma}}) + \frac{1}{2}(\mathbf{x} - \bar{\mathbf{x}})^T\nabla_{\mathbf{xx}}\Xi(\bar{\mathbf{x}},\bar{\bm{\sigma}})(\mathbf{x} - \bar{\mathbf{x}}) \nonumber \\
&=& \frac{1}{2} (\mathbf{x} - \bar{\mathbf{x}})^T G(\bar{\bm{\sigma}}) (\mathbf{x} - \bar{\mathbf{x}}).
\end{eqnarray}
If $G(\bar{\bm{\sigma}})\succeq 0$, it is easy to find that $\bar{\mathbf{x}}$ is the global minimizer of $(\mathcal{P})$, where, $\bar{\mathbf{x}}$ is contained in the canonical equilibrium equation
\begin{eqnarray}
G(\bar{\bm\sigma}) \bar{\mathbf{x}} = F(\bar{\bm\sigma}).
\end{eqnarray}
If $G(\bar{\bm{\sigma}})$ is nonsingular, we have $\bar{\mathbf{x}} = G^{-1}(\bar{\bm\sigma})F(\bar{\bm\sigma})$.
\section{Implementation techniques}
To solve the optimization problem of $(\mathcal{P}^d)$, we firstly rewrite it to the following form:
\begin{eqnarray}
\min && \frac{1}{2}t + \bm\sigma^T \mathbf{d}  \nonumber \\
\mathrm{subject~to:} && t \geq F^T(\bm{\sigma})G^{-1}(\bm{\sigma})F(\bm{\sigma}) \label{eq22} \\
                     && G(\bm{\sigma}) \succeq 0 \label{eq23} \\
                     && \bm{\sigma} \geq 0
\end{eqnarray}
\indent The global solution of $(\mathcal{P}^d)$ is the same to the problem.  Using the Schur complement \cite{Cottle}, we can get the equivalent positive(semi) definite condition to (\ref{eq22}) and (\ref{eq23})
\begin{eqnarray}
\left(
  \begin{array}{cc}
    G(\bm{\sigma}) & F(\bm{\sigma}) \\
    F^T(\bm{\sigma}) & t \\
  \end{array}
\right) \succeq 0
\end{eqnarray}
and then the optimization problem can be expressed as the standard SDP form
\begin{eqnarray}
\min && \frac{1}{2}t + \bm\sigma^T \mathbf{d}  \nonumber \\
\mathrm{subject~to:} &&
\begin{pmatrix}
    G(\bm{\sigma}) & F(\bm{\sigma}) \\
    F^T(\bm{\sigma}) & t
\end{pmatrix} \succeq 0 \\
&& \bm{\sigma} \geq 0
\end{eqnarray}
If $G(\bar{\bm\sigma})\succ 0$, we can get the corresponding global solution to $(\mathcal{P})$ by the canonical duality theory.
In practice, the estimation of $G(\bar{\bm\sigma})$ may exist little inaccuracy due to the perturbed complementary slackness in
primal-dual interior point method and numerical precision. In this study, we use the Cholesky factorization, Condition number and
the smallest Eigenvalue of $G(\bar{\bm\sigma})$ to evaluate the positive definiteness comprehensively. If $G(\bar{\bm\sigma})$ is ill conditioned or det$(G(\bar{\bm\sigma}))=0$, we can add a linear perturbation to the primal objective function and then integrate the canonical dual solutions with the KKT conditions to recover the approximate solution or global solution to primal problem. Details of the techniques are given in the following examples.

\section{Numerical results}
All of the benchmark constrained optimization are from \cite{Liang}, and we keep the number of each problem. In the experiments, we use SeDuMi \cite{Sturn} (a software package which can solve SDP problem) to obtain the canonical dual solutions.
The built-in functions \textit{fsolve} and \textit{fminunc} in MATLAB are also used to solve the simple nonlinear equations and unconstrained optimization problems.\\

\indent \textit{Example 1}: g01
\begin{eqnarray}
\min f(\mathbf{x}) &&= 5 \sum\limits_{i=1}^{4}x_i - 5 \sum\limits_{i=1}^{4}x^2_i - \sum\limits_{i=5}^{13}x_i \nonumber \\
\mathrm{subject~to:} && g_1(\mathbf{x}) = 2x_1 + 2x_2 + x_{10} + x_{11} - 10 \leq 0 \nonumber \\
&& g_2(\mathbf{x}) = 2x_1 + 2x_3 + x_{10} + x_{12} - 10 \leq 0 \nonumber \\
&& g_3(\mathbf{x}) = 2x_2 + 2x_3 + x_{11} + x_{12} - 10 \leq 0 \nonumber \\
&& g_4(\mathbf{x}) = -8x_1 + x_{10} \leq 0 \nonumber \\
&& g_5(\mathbf{x}) = -8x_2 + x_{11} \leq 0 \nonumber \\
&& g_6(\mathbf{x}) = -8x_3 + x_{12} \leq 0 \nonumber \\
&& g_7(\mathbf{x}) = -2x_4 - x_5  + x_{10}\leq 0 \nonumber \\
&& g_8(\mathbf{x}) = -2x_6 - x_7  + x_{11}\leq 0 \nonumber \\
&& g_9(\mathbf{x}) = -2x_8 - x_9  + x_{12}\leq 0 \nonumber
\end{eqnarray}
where the bounds are $0 \leq x_i\leq 1(i = 1, \cdots,9)$, $0 \leq x_i\leq 100(i = 10, 11, 12)$ and
$0 \leq x_{13}\leq 1$.\\
\indent Solving the canonical dual problem, we can obtain $\bar{\bm\sigma}=$
\begin{eqnarray}
\left(
  \begin{array}{c|c|c|c|c|c|c|c|c|c|c}
    \sigma_1 & \sigma_2 & \sigma_3 & \sigma_4 & \sigma_5 & \sigma_6 & \sigma_7 & \sigma_8 & \sigma_9 & \sigma_{10} & \sigma_{11} \\
    0.0000 & 0.0000 & 0.0000 & 0.0000 & 0.0000 & 0.0000 & 1.0000 & 1.0000 & 1.0000 & 5.0000 & 5.0000 \\
    \sigma_{12} & \sigma_{13} & \sigma_{14} & \sigma_{15} & \sigma_{16} & \sigma_{17} & \sigma_{18} & \sigma_{19} & \sigma_{20} & \sigma_{21} & \sigma_{22} \\
    5.0000 & 7.0001 & 2.0001 & 3.0001 & 2.0001 & 3.0001 & 2.0001 & -0.0000 & -0.0000 & -0.0000 & 1.0001 \\
  \end{array}
\right)\nonumber
\end{eqnarray}
\indent In this case, $G(\bar{\bm\sigma}) \succeq 0$ but singular, satisfying the global optimality condition. By the KKT condition, we can find that $g_7,g_8, g_9$, bounds of $x_1, \cdots, x_9$, and $x_{13}$ are active, so we can first get
\begin{eqnarray}
\left(
  \begin{array}{c|c|c|c|c|c|c|c|c|c|c|c|c}
    x_1 & x_2 & x_3 & x_4 & x_5 & x_6 & x_7 & x_8 & x_9 & x_{10} & x_{11} & x_{12} & x_{13} \\
    1 & 1 & 1 & 1 & 1 & 1 & 1 & 1 & 1 & -  & -  & -  & 1  \\
  \end{array}
\right)\nonumber
\end{eqnarray}
where, $``-"$ means undetermined. Considering that constraints $g_7,g_8, g_9$ are active, solving the corresponding linear equations, we can easily get $x_{10} = 3, x_{11} = 3, x_{12} = 3$. Finally, the global solution to g01 is ${\mathbf{x}}^{*}=$
\begin{eqnarray}
\left(
  \begin{array}{c|c|c|c|c|c|c|c|c|c|c|c|c}
    x_1 & x_2 & x_3 & x_4 & x_5 & x_6 & x_7 & x_8 & x_9 & x_{10} & x_{11} & x_{12} & x_{13} \\
    1 & 1 & 1 & 1 & 1 & 1 & 1 & 1 & 1 & 3  & 3  & 3  & 1  \\
  \end{array}
\right)\nonumber
\end{eqnarray}
and $f(\mathbf{x}^{*}) = -15$.\\
\\
\indent \textit{Example 2}: g04
\begin{eqnarray}
\min f(\mathbf{x}) &&= 5.3578547x^2_3 + 0.8356891x_1 x_5 + 37.293239x_1 - 40792.141 \nonumber \\
\mathrm{subject~to:} && g_1(\mathbf{x}) = 85.334407 + 0.0056858x_2x_5 + 0.0006262x_1x_4 - 0.0022053x_3x_5 -92 \leq 0 \nonumber \\
&& g_2(\mathbf{x}) = -85.334407 - 0.0056858x_2x_5 - 0.0006262x_1x_4 + 0.0022053x_3x_5 \leq 0 \nonumber \\
&& g_3(\mathbf{x}) = 80.51249 + 0.0071317x_2x_5 + 0.0029955x_1x_2 + 0.0021813x^2_3 - 110 \leq 0 \nonumber \\
&& g_4(\mathbf{x}) = -80.51249 - 0.0071317x_2x_5 - 0.0029955x_1x_2 - 0.0021813x^2_3 + 90 \leq 0 \nonumber \\
&& g_5(\mathbf{x}) = 9.30096 + 0.0047026x_3x_5 + 0.0012547x_1x_3 + 0.0019085x_3x_4 -25 \leq 0 \nonumber \\
&& g_6(\mathbf{x}) = -9.30096 - 0.0047026x_3x_5 - 0.0012547x_1x_3 - 0.0019085x_3x_4 + 20 \leq 0 \nonumber
\end{eqnarray}
where $78 \leq x_1 \leq 102$, $33 \leq x_2 \leq 45$ and $27 \leq x_i \leq 45(i=3,4,5)$.\\
\indent Solving the canonical dual problem, we can obtain $\bar{\bm\sigma}=$
\begin{eqnarray}
\left(
  \begin{array}{c|c|c|c|c|c|c|c|c|c|c}
    \sigma_1 & \sigma_2 & \sigma_3 & \sigma_4 & \sigma_5 & \sigma_6 & \sigma_7 & \sigma_8 & \sigma_9 & \sigma_{10} & \sigma_{11} \\
    336.8388 & 0.0000 & 0.0001 & 0.0002 & 0.0003 & 798.2826 &  2.0310 & 6.1233 & 0.0001 & 1.6054 & 1.1849 \\
  \end{array}
\right)\nonumber
\end{eqnarray}
\indent In this case, $G(\bar{\bm\sigma}) \succ 0$ and $\mathrm{cond}(G(\bar{\bm\sigma})) = 9.7330e5$, satisfying the global optimality condition, so we can get
$\bar{\mathbf{x}} = $
\begin{eqnarray}
\left(
  \begin{array}{c|c|c|c|c}
    x_1 & x_2 & x_3 & x_4 & x_5 \\
    77.9452 &33.0179 &  29.7345  & 44.9884 &  38.2523\\
  \end{array}
\right)\nonumber
\end{eqnarray}
Noting that the condition number is large, according to the KKT condition, we can first get
\begin{eqnarray}
\left(
  \begin{array}{c|c|c|c|c}
    x_1 & x_2 & x_3 & x_4 & x_5 \\
    78 &33 &  -  & 45 &  -\\
  \end{array}
\right)\nonumber
\end{eqnarray}
Considering that constraints $g_1,g_6$ are active, solving the corresponding linear equations, we can easily get $x_3 = 29.995256025681599, x_5 = 36.775812905788207$. Finally, the global solution to g04 is ${\mathbf{x}}^{*}=$
\begin{eqnarray}
\left(
  \begin{array}{c|c|c|c|c}
    x_1 & x_2 & x_3 & x_4 & x_5 \\
    78 &33 &  29.995256025681599  & 45 &  36.775812905788207\\
  \end{array}
\right)\nonumber
\end{eqnarray}
and $f(\mathbf{x}^{*}) = -3.0666e4$.\\

\noindent \textbf{Remark 1} We use the inverse of $G(\bar{\bm\sigma})$ because only its smallest eigenvalues approximates to zero although its condition number is large. As a matter of fact, the solution $\bar{\mathbf{x}}$ causes only little infeasibility of the first constraint. By integrating the canonical dual solutions and the KKT conditions, we claim that $x_1, x_2$ and $x_4$ are determined in the first stage.
\\

\indent \textit{Example 3}: g07
\begin{eqnarray}
\min f(\mathbf{x}) &&= x^2_1 + x^2_2 + x_1x_2 - 14x_1 - 16x_2 + (x_3 - 10)^2 + 4(x_4 -5)^2 + (x_5 - 3)^2 \nonumber \\
&& ~~~~~~~+2(x_6 -1)^2 + 5x^2_7 + 7(x_8 - 11)^2 + 2(x_9 - 10)^2 + (x_{10} -7)^2 + 45 \nonumber \\
\mathrm{subject~to:} && g_1(\mathbf{x}) = -105 + 4x_1 + 5x_2 - 3x_7 + 9x_8 \leq 0 \nonumber \\
&& g_2(\mathbf{x}) = 10x_1 -8x_2 - 17x_7 + 2x_8 \leq 0 \nonumber \\
&& g_3(\mathbf{x}) = -8x_1 + 2x_2 + 5x_9 - 2x_{10} - 12 \leq 0 \nonumber \\
&& g_4(\mathbf{x}) =  3(x_1 -2)^2 + 4(x_2 -3)^2 + 2x_3^2 - 7x_4 - 120\leq 0 \nonumber \\
&& g_5(\mathbf{x}) =  5x_1^2 + 8x_2 + (x_3 -6)^2 - 2x_4 - 40\leq 0 \nonumber \\
&& g_6(\mathbf{x}) =  x_1^2 + 2(x_2 -2)^2 - 2x_1 x_2 + 14x_5 - 6x_6\leq 0 \nonumber \\
&& g_7(\mathbf{x}) =  0.5(x_1 - 8)^2 + 2(x_2 - 4)^2 + 3 x_5^2 - x_6 -30 \leq 0 \nonumber \\
&& g_8(\mathbf{x}) =  -3x_1 + 6x_2 + 12(x_9 - 8)^2 - 7 x_{10}\leq 0 \nonumber
\end{eqnarray}
where $-10 \leq x_i \leq 10 (i = 1, \cdots, 10)$.\\
\indent Solving the canonical dual problem, we can obtain $\bar{\bm\sigma}=$
\begin{eqnarray}
\left(
  \begin{array}{c|c|c|c|c|c|c|c|c|}
    \sigma_1 & \sigma_2 & \sigma_3 & \sigma_4 & \sigma_5 & \sigma_6 & \sigma_7 & \sigma_8 & \sigma_9 \\
    1.7168 & 0.4746 & 1.3760 & 0.0205 & 0.3120 & 0.2871 &  0.0000 & 0.0000 & 0.0000 \\
    \sigma_{10} & \sigma_{11} & \sigma_{12} & \sigma_{13} & \sigma_{14} & \sigma_{15} & \sigma_{16} & \sigma_{17} & \sigma_{18} \\
    0.0000 & 0.0000 & 0.0000 & 0.0000 & 0.0000 & 0.0000 &  0.0000 & 0.0000 & 0.0000\\
  \end{array}
\right)\nonumber
\end{eqnarray}
\indent In this case, $G(\bar{\bm\sigma}) \succ 0$ and $\mathrm{cond}(G(\bar{\bm\sigma}))= 7.0000$, satisfying the global optimality condition, so we can get
$\mathbf{x}^{*} = $
\begin{eqnarray}
\left(
  \begin{array}{c|c|c|c|c|c|c|c|c|c}
    x_1 & x_2 & x_3 & x_4 & x_5 & x_6 & x_7 & x_8 & x_9 & x_{10}\\
    2.1721 & 2.3636 &   8.7746  &  5.0959 &   0.9903 &   1.4307 &   1.3218  &  9.8286 &   8.2800 &   8.3760\\
  \end{array}
\right)\nonumber
\end{eqnarray}
and $f(\mathbf{x}^{*}) = 24.3111$.
Note that there exists little infeasibility due to numerical precision.\\
\indent \textit{Example 4}: g10
\begin{eqnarray}
\min f(\mathbf{x}) &&= x_1 + x_2 + x_3 \nonumber \\
\mathrm{subject~to:} && g_1(\mathbf{x}) = -1 + 0.0025(x_4 + x_6) \leq 0 \nonumber \\
&& g_2(\mathbf{x}) = -1 + 0.0025(x_5 + x_7 - x_4) \leq 0 \nonumber \\
&& g_3(\mathbf{x}) = -1 + 0.01(x_8 - x_5) \leq 0 \nonumber \\
&& g_4(\mathbf{x}) =  -x_1 x_6 + 833.33252x_4 + 100x_1 - 83333.333 \leq 0 \nonumber \\
&& g_5(\mathbf{x}) =  -x_2 x_7 + 1250x_5 + x_2 x_4 - 1250x_4 \leq 0 \nonumber \\
&& g_6(\mathbf{x}) =  -x_3 x_8 + 1250000 + x_3 x_5 - 2500x_5 \leq 0 \nonumber
\end{eqnarray}
where $100 \leq x_1 \leq 10000$, $1000 \leq x_i \leq 10000 (i = 2, 3)$ and $10 \leq x_i \leq 1000 (i = 4, \cdots, 8)$\\
\indent Solving the canonical dual problem, we can obtain $\bar{\bm\sigma}=$
\begin{eqnarray}
\left(
  \begin{array}{c|c|c|c|c|c|c|}
    \sigma_1 & \sigma_2 & \sigma_3 & \sigma_4 & \sigma_5 & \sigma_6 & \sigma_7\\
    9.2834  &28.9205  & 5.8893   & 0.0001 &   0.0001  &  0.0001  &  0.0001\\
    \sigma_{8} & \sigma_{9} & \sigma_{10} & \sigma_{11} & \sigma_{12} & \sigma_{13} & \sigma_{14}\\
   0.0001  &  0.0001 &   0.0000 &  0.0000  &  0.0000 &   0.0000  &  0.0000
  \end{array}
\right)\nonumber
\end{eqnarray}
\indent In this case, $G(\bar{\bm\sigma}) \succ 0$ and $\mathrm{cond}(G(\bar{\bm\sigma})) = 749.4514$, satisfying the global optimality condition. However, the $\mathrm{max}(\mathrm{eig}(G(\bar{\bm\sigma}))) =  2.5743e-4$, which is too small, so we cannot use the inverse of
$G(\bar{\bm\sigma})$ directly. By the KKT condition, we can find that constraints $g_1, g_2, g_3$ are active, and all of the box constraints are inactive. That is to say, the problem is equivalent to a linear programming problem with linear constraints, which indicates that $g_4, g_5, g_6$ must be active.  Fixing $x_4, x_5$, we have
\begin{equation*}
\left\{ \begin{aligned}
x_1 & = \frac{83333.333 - 833.33252x_4}{x_4 - 300}\\
x_2 & =  \frac{1250x_4 - 1250x_5}{x_5 - 400}\\
x_3 & = 12500 - 25x_5\\
x_6 & = 400 - x_4\\
x_7 & = 400 + x_4 - x_5\\
x_8 & = 100 + x_5
\end{aligned} \right.
\end{equation*}
As a result, we can reduce the problem to
\begin{eqnarray*}
\min f(\mathbf{x}) &&=  \frac{83333.333 - 833.33252x_4}{x_4 - 300} + \frac{1250x_4 - 1250x_5}{x_5 - 400} + 12500 - 25x_5
\end{eqnarray*}
Taking the box constraints of $x_1, \cdots, x_8$ into consideration, when using $(100, 200)$ as an initial point for the unconstrained optimization problem with two variables, it is easy to get the only minimum $x_4 = 182.0176995811199$ and $x_5 = 295.6011732779338$. Utilizing the equations obtained by the complementary slackness, finally, we have $x_1 = 579.3066844253549$
$x_2 = 1359.970668051655$, $x_3 = 5109.970668051655$, $x_6 = 217.9823004188801$, $x_7 = 286.4165263031861$, $x_8 = 395.6011732779338$
and $f(\mathbf{x}^{*}) = 7049.248020528666$.\\

\noindent \textbf{Remark 2} We don't use the inverse of $G(\bar{\bm\sigma})$ directly because all of its eigenvalues are approximately zeros. And the reason why we still use the canonical dual solutions as useful heuristics is that the $G(\bar{\bm\sigma})$ is slightly positive definite due to the perturbed complementary slackness caused by the SeDuMi. Since all of the box constraints are inactive and the objective function is linear, it is not difficult to imagine that all of the constraints must be active. Note that the constraints of $x_4$ and $x_5$ are changed when solving the unconstrained optimization problem since constraints of $x_1, x_2, x_3$ and $x_6, x_7, x_8$ must be satisfied.
\\

\indent \textit{Example 5}: g18
\begin{eqnarray}
\min f(\mathbf{x}) &&= -0.5(x_1 x_4 - x_2 x_3 + x_3 x_9 - x_5 x_9 + x_5 x_8 - x_6 x_7) \nonumber \\
\mathrm{subject~to:} && g_1(\mathbf{x}) = x_3^2 + x_4^2 - 1 \leq 0 \nonumber \\
&& g_2(\mathbf{x}) = x_9^2 -1 \leq 0 \nonumber \\
&& g_3(\mathbf{x}) = x_5^2 + x_6^2 -1 \leq 0 \nonumber \\
&& g_4(\mathbf{x}) = x_1^2 + (x_2 - x_9)^2 -1 \leq 0 \nonumber \\
&& g_5(\mathbf{x}) =  (x_1 - x_5)^2 + (x_2 - x_6)^2 -1 \leq 0 \nonumber \\
&& g_6(\mathbf{x}) =  (x_1 - x_7)^2 + (x_2 - x_8)^2 - 1 \leq 0 \nonumber \\
&& g_7(\mathbf{x}) =  (x_3 - x_5)^2 + (x_4 - x_6)^2 - 1 \leq 0 \nonumber \\
&& g_8(\mathbf{x}) =  (x_3 - x_7)^2 + (x_4 - x_8)^2 - 1 \leq 0 \nonumber \\
&& g_9(\mathbf{x}) =  x_7^2 + (x_8 - x_9)^2 - 1 \leq 0 \nonumber \\
&& g_{10}(\mathbf{x}) =  x_2x_3 - x_1x_4  \leq 0 \nonumber \\
&& g_{11}(\mathbf{x}) =  -x_3 x_9  \leq 0 \nonumber \\
&& g_{12}(\mathbf{x}) =  x_5 x_9  \leq 0 \nonumber \\
&& g_{13}(\mathbf{x}) =  x_6 x_7 - x_5x_8  \leq 0 \nonumber
\end{eqnarray}
where $-10 \leq x_1 \leq 10, (i = 1, \cdots, 8)$ and $0 \leq  x_9 \leq 20$.\\
\indent Solving the canonical dual problem, we can obtain $\bar{\bm\sigma}=$
\begin{eqnarray}
\left(
  \begin{array}{c|c|c|c|c|c|c|c|c|c|c}
    \sigma_1 & \sigma_2 & \sigma_3 & \sigma_4 & \sigma_5 & \sigma_6 & \sigma_7 & \sigma_8 & \sigma_9 & \sigma_{10} & \sigma_{11} \\
    0.1444  &  0.0000  &  0.1444   & 0.1445  &  0.0000  &  0.1442  &  0.1441 &   0.0000 &   0.1445 &   0.0000  &  0.0000 \\
    \sigma_{12} & \sigma_{13} & \sigma_{14} & \sigma_{15} & \sigma_{16} & \sigma_{17} & \sigma_{18} & \sigma_{19} & \sigma_{20} & \sigma_{21} & \sigma_{22} \\
     0.0000 &   0.0000 &  -0.0000 &  -0.0000 &  -0.0000&   -0.0000 &  -0.0000  & -0.0000 &  -0.0000  &  -0.0000  &  0.0000 \\
  \end{array}
\right)\nonumber
\end{eqnarray}
\indent In this case, $G(\bar{\bm\sigma}) \succ 0$ and $\mathrm{cond}(G(\bar{\bm\sigma})) = 7.1887e7$, satisfying the global optimality condition. However, the condition number is large. Taking the KKT conditions into account, we can conclude that constraints $g_1$, $g_3$, $g_4$, $g_6$, $g_7$, $g_9$ are active since the corresponding $\sigma_1$, $\sigma_3$, $\sigma_4$, $\sigma_6$, $\sigma_7$,  $\sigma_9$ are not zeros. But it becomes still difficult to solve the nonlinear equations.
Considering that several eigenvalues of $G(\bar{\bm\sigma})$ are zeros and there exists no linear term in the objective function, and in this situation, we add a small linear perturbation $0.05 (x_1 + \cdots, x_9) $ to the primal objective function. Solving the perturbed canonical dual problem, we get $G(\bar{\bm\sigma}) \succ 0$ and $\mathrm{cond}(G(\bar{\bm\sigma})) = 1.4592e3$ and the smallest eigenvalue of $G(\bar{\bm\sigma})$ is 0.0021. Therefore,
we can get $\bar{\mathbf{x}}$ =
\begin{eqnarray}
\left(
  \begin{array}{c|c|c|c|c|c|c|c|c|c}
    x_1 & x_2 & x_3 & x_4 & x_5 & x_6 & x_7 & x_8 & x_9\\
     -0.9660  & -0.2585 &  -0.2587 &  -0.9660 &  -0.9661  & -0.2588 &  -0.2589 &  -0.9657 &    0.0005
  \end{array}
\right)\nonumber
\end{eqnarray}
and $f(\bar{\mathbf{x}}) = -0.8663$. Note that there exists little infeasibility due to numerical precision.\\

\noindent \textbf{Remark 3}  The solution we get is quite different from the best known solution. According to the canonical duality theory, the global solution to this problem is not unique. The linear perturbation technique can only help to find one
of the global solutions. As a matter of fact, the following solutions $\bar{\mathbf{x}}$ =
\begin{eqnarray}
\left(
  \begin{array}{c|c|c|c|c|c|c|c|c|c}
    x_1 & x_2 & x_3 & x_4 & x_5 & x_6 & x_7 & x_8 & x_9\\
    0.0450 &  -0.0387  &   0.8663 &  -0.4999  &  0.0004  & -1.0001 &   0.8878 &   0.5000  &  0.9604\\
  \end{array}
\right),\nonumber
\end{eqnarray}
\begin{eqnarray}
\left(
  \begin{array}{c|c|c|c|c|c|c|c|c|c}
    x_1 & x_2 & x_3 & x_4 & x_5 & x_6 & x_7 & x_8 & x_9\\
     0.0689  &   -0.9972  &    0.9088  &   -0.4179  &    0.0920   &  -0.9959   &   0.8986  &   -0.4388  &    0.0009\\
  \end{array}
\right),\nonumber
\end{eqnarray}
and
\begin{eqnarray}
\left(
  \begin{array}{c|c|c|c|c|c|c|c|c|c}
    x_1 & x_2 & x_3 & x_4 & x_5 & x_6 & x_7 & x_8 & x_9\\
    0.6888 &  -0.7257  &  0.9693  &  0.2454 &   0.6973 &  -0.7173  &  0.9726  &  0.2332 &  -0.0006\\
  \end{array}
\right)\nonumber
\end{eqnarray}
can all be considered as approximate solutions, which are obtained by the proposed techniques.
\section{Conclusion}
We have applied the canonical duality theory to solve a class of CEC benchmark constrained optimization problems.
Experimental results show that some of the examples can be solved directly, some of them can be solved by integrating the canonical dual solutions and KKT conditions, and some can be solved approximately by adding a small linear perturbation term.




\end{document}